 \documentclass[reqno]{amsart} \usepackage{amscd}
 \usepackage{epsf}
\usepackage{epsfig}
 \usepackage{graphicx}
 \usepackage{tikz} 
 \usepackage{amssymb}
 \usepackage{tkz-euclide}
 \usepackage{verbatim}
 \usepackage[english]{babel}
 \usepackage{tabu}
 \usepackage{amsmath}
 \newtheorem{theorem}{Theorem}[section]

 \newtheorem{lemma}[theorem]{Lemma}

 \theoremstyle{remark} 
 \numberwithin{equation}{section}

 \providecommand \hf{\hspace*{0.5cm}}
 
 \theoremstyle{plain}
 \newtheorem{thm}{Theorem}[section]

 \newtheorem*{que}{Question}
 
 \newtheorem{rem}[thm]{Remark}
 \newtheorem*{mresult}{Main Result}

 \theoremstyle{definition}

 \usepackage{hyperref}
 
\providecommand \hf{\hspace*{0.5cm}}

\begin{document}
\title{Rational cuspidal curves in a moving family of $\mathbb{P}^2$}

\author[Ritwik Mukherjee]{Ritwik Mukherjee}
\address{School of Mathematics, National Institute of Science Education and Research, Bhubaneswar, HBNI, Odisha 752050, India}
\email{ritwikm@niser.ac.in}
 
\author[Rahul Kumar Singh]{Rahul Kumar Singh}
\address{Department of Mathematics, Indian Institute of Technology Patna, Bihta, Patna-801106, India}
\email{rahulks@iitp.ac.in}
\subjclass[2010]{14N35, 14J45}

\begin{abstract}
In this paper we obtain a formula for  
the number of rational degree $d$ curves in $\mathbb{P}^3$ having a cusp, whose image lies in a $\mathbb{P}^2$ and that passes through $r$ 
lines and $s$ points (where $r + 2s = 3d+1$). 
This problem can be viewed as a family version of the classical question of counting rational cuspidal curves in $\mathbb{P}^2$, 
which has been studied earlier by Z. Ran (\cite{RanBB2}), R.~Pandharipande (\cite{RP}) and A. Zinger (\cite{g2p2and3}).
We obtain this number by computing the Euler class of a relevant bundle and then finding out the corresponding 
degenerate contribution to the Euler class. The method we use is closely based on the method followed 
by A. Zinger (\cite{g2p2and3}) and I. Biswas, S. D'Mello, R. Mukherjee and V. Pingali (\cite{BSRV}). 
We also verify that our answer for 
the characteristic numbers of rational cuspidal planar cubics and quartics is 
consistent with the answer obtained by N. Das and the first author (\cite{ND}), where they 
compute the characteristic 
number of $\delta$-nodal planar curves in $\mathbb{P}^3$ with one cusp
(for $\delta \leq 2$).
\end{abstract}

\maketitle
\tableofcontents

\section{Introduction}
\noindent A classical question in enumerative algebraic geometry is: 
\begin{que}
\label{q1}
What is $N_d$, the number of rational (genus zero) degree $d$ curves in $\mathbb{P}^2$ that pass through $3d-1$ generic points?  
\end{que}
\hf \hf Although the computation of $N_d$ 
is a classical question, a complete solution to the above problem 
was unknown until the early
$90^{' \textnormal{s}}$ when Ruan--Tian (\cite{RT}) and Kontsevich--Manin (\cite{KoMa})
obtained a formula for $N_d$. Generalization of this question to enumerate rational curves with higher singularities 
(such as cusps, tacnodes and higher order cusps) have been studied by Z. Ran (\cite{RanBB2}), R.~Pandharipande (\cite{RP}) and A. Zinger (\cite{g2p2and3}, 
\cite{genus_three_zinger} and \cite{g0pr}). These results have also been generalized to other surfaces (such as $\mathbb{P}^1 \times \mathbb{P}^1$) 
by J. Kock (\cite{kock}) and more generally for del-Pezzo surfaces by I. Biswas, S. D'Mello, R. Mukherjee and V. Pingali (\cite{BSRV}).
The problem of enumerating elliptic cuspidal curves has been solved by Z. Ran (\cite{RanBB2}), and more recently a solution to this question 
in any genus has been obtained 
by Y. Ganor and E. Shustin 
(\cite{GS}) using methods from Tropical Geometry. \\ 
\hf \hf A natural generalization of problems in enumerative geometry (where one studies curves inside some fixed ambient surface such as $\mathbb{P}^2$) 
is to consider a family version of the same problem. This generalization is considered by S. Kleiman and R. Piene (\cite{KP2}) 
and more recently by T.~Laarakker (\cite{TL})
where they study the 
enumerative geometry of nodal curves in a moving family of 
surfaces. \\
\hf \hf Motivated by this generalization, A. Paul and the authors of this paper 
studied a family version of computing $N_d$ in (\cite{ritwik}); 
there the authors find a formula for the characteristic number of rational planar curves in $\mathbb{P}^3$ 
(i.e. curves in $\mathbb{P}^3$ that lie inside a $\mathbb{P}^2$). 
In this paper we build up on the results of (\cite{ritwik}) to find the characteristic number of rational planar curves in $\mathbb{P}^3$ having a cusp. 
The main result of this paper is as follows: 
\begin{mresult}
Let $C_{d}^{\mathbb{P}^3, \mathrm{Planar}}(r,s)$ be the number of genus zero, degree $d$ curves in $\mathbb{P}^3$ having a cusp, 
whose image  lies in a $\mathbb{P}^2$, intersecting 
$r$ generic lines and $s$ generic points (where $r  + 2s = 3d+1$). We have 
a recursive formula to compute $C_{d}^{\mathbb{P}^3, \mathrm{Planar}}(r,s)$. 
\end{mresult}
\noindent We have written a mathematica program to implement our formula; 
the program is available on our web page
\[ \textnormal{\url{https://www.sites.google.com/site/ritwik371/home}}. \]
\hf \hf In section \ref{low_deg_chck}, we subject our formula to several low degree checks; in particular, we 
verify that our numbers are logically consistent with those obtained by N. Das and the first author (\cite{ND}). \\ 
\hf \hf Let us now give a brief overview of the method we use in this paper; we closely adapt the method 
applied by A. Zinger (\cite{g2p2and3}) and  
I. Biswas, S. D'Mello, R. Mukherjee and V. Pingali (\cite{BSRV}). 
We express our enumerative number as the number of zeros 
of a section of an appropriate vector bundle (restricted to an open dense set of an appropriate moduli space). 
As is usually the case, the Euler class of this vector bundle is our desired enumerative number, \textit{plus} 
an extra boundary contribution. In (\cite{g2p2and3}) and (\cite{BSRV}) 
the method of ``dynamic intersections'' (cf.~Chapter 11 in \cite{F}) is used 
to compute the degenerate contribution to the relevant Euler class.
We argue in section \ref{Degenerate_Euler} how the multiplicity computation in 
(\cite{g2p2and3}) and (\cite{BSRV}) implies the multiplicity of the degenerate locus 
that occurs in our case. 
Finally, computation of the Euler class involves the intersection of tautological classes on the moduli space of planar degree $d$ curves; this 
in turn involves the characteristic number of rational planar curves in $\mathbb{P}^3$, for which we use the result of our paper \cite{ritwik}. 
Hence, we can compute both the Euler class and the degenerate contribution, which gives us our desired number $C_{d}^{\mathbb{P}^3, \mathrm{Planar}}(r,s)$.  


\section{Notation}
\label{notation}
Let us define a \textbf{planar} curve in $\mathbb{P}^3$ to be a curve, whose image lies inside a $\mathbb{P}^2$. 
We will now develop some notation to describe the space of planar curves of a given degree $d$. \\ 
\hf \hf Let us denote the dual of $\mathbb{P}^3$ by 
$\widehat{\mathbb{P}}^3$; this is the space of $\mathbb{P}^2$ inside $\mathbb{P}^3$. 
An element of $\widehat{\mathbb{P}}^3$ can be thought of as 
a non zero linear functional $\eta: \mathbb{C}^4 \longrightarrow \mathbb{C}$ upto scaling (i.e., it is the projectivization of the dual of $\mathbb{C}^4$).  
Given such an $\eta$, we  define the projectivization of its zero set as $\mathbb{P}^2_{\eta}$. In other words, 
\begin{align*}
\mathbb{P}^2_{\eta} &:= \mathbb{P}(\eta^{-1}(0)). 
\end{align*}
Note that this $\mathbb{P}^2_{\eta}$ is a subset of $\mathbb{P}^3$. Next, 
we define the moduli space of planar degree $d$ curves into $\mathbb{P}^3$ 
as a fibre bundle over $\widehat{\mathbb{P}}^3$. More precisely, 
we define 
\[ \pi: \overline{\mathcal{M}}_{0,1}^{\mathrm{Planar}}(\mathbb{P}^3,d) \longrightarrow \widehat{\mathbb{P}}^3 \] 
to be the fiber bundle, such that 
\[ \pi^{-1}([\eta]) := \overline{\mathcal{M}}_{0,1}(\mathbb{P}^2_{\eta}, d). \]
Here we are using the standard notation to denote $\mathcal{M}_{0,k}(X, \beta)$ to be the moduli space of genus zero stable maps, 
representing the class $\beta \in H_2(X, \mathbb{Z})$ and $\overline{\mathcal{M}}_{0,k}(X, \beta)$ to be its stable map compactification. 
Since the dimension of a fiber bundle is the dimension of the base, plus the dimension of the fiber, we conclude that the dimension of 
$\overline{\mathcal{M}}_{0,k}^{\mathrm{Planar}}(\mathbb{P}^3,d)$ is $3d+2+k$.\\
\hf \hf Next, we note that there is a natural forgetful map 
\begin{align*}
\pi_{F}: \overline{\mathcal{M}}_{0,1}^{\mathrm{Planar}}(\mathbb{P}^3,d) & \longrightarrow \overline{\mathcal{M}}_{0,1}(\mathbb{P}^3, d)  
\end{align*}
where one forgets the plane $\mathbb{P}^2_{\eta}$ 
and simply thinks about the stable map to $\mathbb{P}^3$. 
When $d \geq 2$, the map $\pi_{F}$ is injective when restricted to the open dense subspace of non multiply covered curves 
(from a smooth domain). 
This is because every planar degree $d$ degree curve lies in a unique 
plane, when $d \geq 2$. When $d=1$, this map is not injective since a line is not contained in a unique plane.  
Infact we note that the space of  lines is $4$ 
dimensional, while the dimension of $\overline{\mathcal{M}}_{0,0}^{\mathrm{Planar}}(\mathbb{P}^3,1)$ is $5$.\\
\hf \hf Let 
\[\mathbb{L} \longrightarrow \overline{\mathcal{M}}_{0,1}(\mathbb{P}^3, d)\] 
denote the universal tangent line bundle over the marked point (the fiber over each point is the tangent space over the given marked point). 
This line bundle will pullback to a  
line bundle over $\overline{\mathcal{M}}_{0,1}^{\mathrm{Planar}}(\mathbb{P}^3,d)$ via the map $\pi_{F}$; we will denote it by the same 
symbol $\mathbb{L}$ (we will in general avoid writing the pullback symbol $\pi_F^*$ if there is no cause of confusion). \\ 
\hf \hf Let us now define a few cycles in $\overline{\mathcal{M}}_{0,1}^{\mathrm{Planar}}(\mathbb{P}^3,d)$. 
Let 
$\mathcal{H}_L$ and $\mathcal{H}_p$  denote the classes  of the cycles in $\overline{\mathcal{M}}_{0,0}^{\mathrm{Planar}}(\mathbb{P}^3,d)$ that corresponds 
to the subspace of curves passing through a generic line and a point respectively. 
We will denote their pullbacks (via the forgetful map that forgets the marked point) to 
$\overline{\mathcal{M}}_{0,0}^{\mathrm{Planar}}(\mathbb{P}^3,d)$ by the same symbol  
$\mathcal{H}_L$ and $\mathcal{H}_p$. 
We will also denote  $H$ and $a$ 
to be the standard generators of $H^*(\mathbb{P}^3; \mathbb{Z})$ and $H^*(\widehat{\mathbb{P}}^3; \mathbb{Z})$ respectively. 
As $\pi$ is a projection map from $\overline{\mathcal{M}}_{0,1}^{\mathrm{Planar}}(\mathbb{P}^3,d)$ to $\widehat{\mathbb{P}}^3$; 
we denote the pullback $\pi^*a$ by the same symbol $a$. Finally, there is an evaluation map from 
\begin{align*}
\textnormal{ev}:\overline{\mathcal{M}}_{0,1}^{\mathrm{Planar}}(\mathbb{P}^3,d)& \longrightarrow \mathbb{P}^3. 
\end{align*}
We will denote the pullback of $H$ via this map to be $\textnormal{ev}^*H$. This is the one case where we explicitly write the 
pullback map (since it will avoid confusion; in the remaining cases omitting to write down the pullback map causes no confusion).\\  
\hf \hf We will now define a few numbers by intersecting cycles on $\overline{\mathcal{M}}_{0,k}^{\mathrm{Planar}}(\mathbb{P}^3,d)$.
We will use the convention that 
\begin{align}
\langle \alpha, ~[M] \rangle & =0  \qquad \textnormal{if} \qquad \textnormal{deg}(\alpha) \neq \textnormal{dim}(M). \label{convention_deg}
\end{align}
We now define 
\begin{align}
N_d^{\mathbb{P}^3,\mathrm{Planar}}(r,s,\theta) &:= \langle a^{\theta}, ~~\overline{\mathcal{M}}_{0,0}^{\mathrm{Planar}}(\mathbb{P}^3,d) 
\cap \mathcal{H}_L^r \cap \mathcal{H}_p^s \rangle, \nonumber \\   
\Phi_d(i,j, r, s,  \theta) &:= 
\langle c_1(\mathbb{L}^*)^i(\textnormal{ev}^*H)^j,\; \overline{\mathcal{M}}_{0,1}^{\mathrm{Planar}}(\mathbb{P}^3,d)\cap 
\mathcal{H}_L^r \cap \mathcal{H}_p^s\cap a^{\theta}\rangle. \nonumber 
\end{align}
We note that using the results of our paper (\cite{ritwik}), the numbers $N_d^{\mathbb{P}^3,\mathrm{Planar}}(r,s,\theta)$ 
are all computable. In section \ref{itc}, a formula is given to compute the relevant $\Phi_d(i,j, r, s, \theta)$ necessary 
to obtain the main result  of this paper. We note that using the convention introduced in equation \eqref{convention_deg} 
\begin{align*}
N_d^{\mathbb{P}^3,\mathrm{Planar}}(r,s,\theta) & = 0 \qquad \textnormal{unless} \qquad r+2s+\theta = 3d+2 \qquad \textnormal{and} \\ 
\Phi_d(i,j, r, s,  \theta) & = 0 \qquad \textnormal{unless} \qquad r+2s+\theta + i + j = 3d+3. 
\end{align*}

\section{Euler class computation}
\label{Euler_class_computation}
\noindent We will now describe the basic method by which we compute 
the characteristic number of planar rational cuspidal curves in $\mathbb{P}^3$. 
We will  express this number as the 
number of zeros of a section of an appropriate bundle restricted to an open dense subspace 
of the moduli space of planar curves in $\mathbb{P}^3$. \\ 
\hf \hf Before we do that, let us make a few abbreviations that 
we will often use. We denote 
\begin{align*}
C_d&:= C_d^{\mathbb{P}^3, \textnormal{Planar}}(r,s), \qquad 
\mathcal{M}:= \mathcal{M}_{0,1}^{\textnormal{Planar}}(\mathbb{P}^3, d) \cap  
\mathcal{H}_L^r \cap \mathcal{H}_p^s \qquad \textnormal{and} \\ 
\overline{\mathcal{M}}& := \overline{\mathcal{M}}_{0,1}^{\textnormal{Planar}}(\mathbb{P}^3, d) \cap  
\mathcal{H}_L^r \cap \mathcal{H}_p^s.
\end{align*}
\noindent Let us now define 
\begin{align*}
\mathcal{S}&:= \{ ([\eta],  q) \in \widehat{\mathbb{P}}^3 \times \mathbb{P}^3: \eta(q)=0\}. 
\end{align*}
An element of $\mathcal{S}$ denotes a plane $\mathbb{P}^2_{\eta}$ in $\mathbb{P}^3$ together with a marked point $q$ that 
lies in the plane. 
We will now define a rank two bundle $W \longrightarrow \mathcal{S}$, where the fibre over each point $([\eta], q)$ is 
$T_q\mathbb{P}^2_{\eta}$. 
We note that over $\mathcal{S}$, we have the following short exact sequence of vector bundles:  
\begin{align*}
0&\longrightarrow W \longrightarrow T \mathbb{P}^3 \longrightarrow \gamma^{*}_{\widehat{\mathbb{P}}^3} \otimes \gamma^{*}_{\mathbb{P}^3} \longrightarrow 0. 
\end{align*}
Here the first map is the inclusion map and the second map is $\nabla \eta|_q$.
Hence, 
\begin{align*}
c(W) c(\gamma^{*}_{\widehat{\mathbb{P}}^3} \otimes \gamma^{*}_{\mathbb{P}^3}) & = c(T\mathbb{P}^3) ~~
\implies c_1(W) = 3H-a \qquad \textnormal{and} \\
\qquad c_2(W) = a^2-2aH + 3H^2.
\end{align*}
Next, we note that $C_d$ is the cardinality of the set 
\begin{align*}
\{ [u, y] \in \mathcal{M}: du|_y = 0\}. 
\end{align*}
The process of taking the derivative of $u$ at the marked point induces a section of the rank two vector bundle 
\begin{align}
\mathbb{L}^* \otimes \textnormal{ev}^*W \longrightarrow \overline{\mathcal{M}}, 
\qquad \textnormal{given by} \qquad [u, y] \longrightarrow du|_y. \label{du_transverse}
\end{align}
This section is transverse to zero, when restricted to $\mathcal{M}$; this is justified 
in section \ref{Degenerate_Euler}. 
We now note that there is a degenerate contribution to the Euler class because the 
section vanishes on the boundary $\overline{\mathcal{M}} - \mathcal{M}$. 
Let us first describe the boundary component on which it will vanish. The boundary component on 
which it will vanish is going to be a map from a wedge of three spheres of degree $d_1$, $0$ and $d_2$, where the marked point lies on the 
degree zero component (which is also called a ghost bubble). 
The cardinality of this set is computed in (\cite{ritwik})  (in the proof of Theorem 3.3); it is given by 
\begin{align}
\mathrm{B} & = \frac{1}{2} \sum_{\substack{d_1+d_2 = d, \\ s_1+s_2 = s, \\ r_1+r_2 = r}} d_1 d_2 B_{d_1, d_2}(r_1,s_1,r_2, s_2, \theta)\binom{r}{r_1} \binom{s}{s_1},  
\label{B_value}
\end{align}
where $B_{d_1, d_2}(r_1,s_1,r_2, s_2, \theta)$ is as defined in equation \eqref{bd1d2_simple}.
We claim that this boundary contributes with a multiplicity of one to the Euler class; this is justified in section \ref{Degenerate_Euler}.\\
\hf \hf It remains to compute the Euler class. 
By the splitting principle, we note that 
\begin{align}
e:= \langle e(\mathbb{L}^* \otimes \textnormal{ev}^*W), \overline{\mathcal{M}}\rangle  & = 
\langle c_1(\mathbb{L}^*)^2 + c_1(\mathbb{L}^*)c_1(\textnormal{ev}^*W) + c_2(\textnormal{ev}^*W), \overline{\mathcal{M}} \rangle \nonumber \\ 
& = \Phi_{d}(2,0,r,s, 0) - \Phi_d(1,0,r ,s ,1) + 3 \Phi_d(1,1 ,r ,s ,0 ) \nonumber \\ 
& + \Phi_d(0,0 ,r ,s ,2)-2 \Phi_d(0,1 ,r ,s ,1) \nonumber \\ 
& + 3 \Phi_d(0,2 ,r ,s ,0). \label{Euler_ITC}
\end{align}
The numbers $\Phi_d(i,j, r, s, \theta)$ that arise in the right hand side of equation \eqref{Euler_ITC} 
can be computed using the results of section \ref{itc} (namely, Lemmas \ref{lemma1}, \ref{l1} and \ref{l2}). 
Since the boundary contributes with a multiplicity of one, we conclude that 
\begin{align}
e &= C_d +  \mathrm{B}. \label{Euler_Cd} 
\end{align}
Using equations \eqref{Euler_Cd}, \eqref{Euler_ITC}, the values of $\Phi_d(i,j, r, s, \theta)$ from the results of section \ref{itc}, 
equation \eqref{B_value} and the values of 
$N_{d}^{\mathbb{P}^3, \mathrm{Planar}}(r,s, \theta)$ from the paper (\cite{ritwik}), we obtain the value of 
$C_d^{\mathbb{P}^3, \textnormal{Planar}}(r,s,\theta)$. This formula has been implemented using mathematica; 
the program is available on request. In section \ref{low_deg_chck}, we present the values of 
$C_d$ for a few values of $d$.

\section{Intersection of Tautological Classes}
\label{itc}
\noindent We  will now give a formula for the relevant $\Phi_d(i,j, r, s, \theta)$ that are necessary to compute the Euler  class. 
We will often refer to $\Phi_d(i,j, r, s, \theta)$ as a level $i$ number. 

\begin{lemma}\label{lemma1}
The level zero numbers $\Phi_d(0,j,r,s,\theta)$ are given by 
\begin{align}
\label{l0}
\Phi_d(0,j,r,s,\theta) & = \begin{cases} 0 & \mbox{if} ~~j=0, \\ 
dN_d^{\mathbb{P}^3, \mathrm{Planar}}(r,s,\theta) & \mbox{if }  ~~j=1,\\ 
N_d^{\mathbb{P}^3, \mathrm{Planar}}(r+1,s,\theta) & \mbox{if } ~~j=2,\\ 
N_d^{\mathbb{P}^3, \mathrm{Planar}}(r,s+1,\theta) &  \mbox{if} ~~ j=3,\\ 
0 & \mbox{if} ~~ j>3.\\ 
\end{cases} 
\end{align}
\end{lemma}
 
\begin{lemma}
\label{l1}
The level one numbers $\Phi_d(1,j,r,s,\theta)$ are given by 
\begin{align}
\Phi_d(1,j,r,s,\theta) & = \begin{cases} -2N_d(r,s,\theta) & \mbox{if} ~~j=0, \\ 
\frac{1}{d^2}\Phi_d(0,1,r+1,s,\theta)-\frac{2}{d}\Phi_d(0,2,r,s,\theta)\; \\
+\frac{1}{d^2}\sum\limits_{\substack{r_1+r_2=r \\ s_1+s_2=s \\ d_1+d_2=d \\ 
d_1,d_2>0}} d_1^2d_2^3\binom{r}{r_1}\binom{s}{s_1} B_{d_1,d_2}(r_1,s_1,r_2,s_2,\theta) & \mbox{if }  ~~j=1,
\end{cases} 
\end{align}
where $B_{d_1,d_2}(r_1,s_1,r_2,s_2,\theta)$ is as defined in equation \eqref{bd1d2_simple}. 
\end{lemma}


\begin{lemma}
\label{l2}
The level two number $\Phi_d(2,0,r,s,\theta)$ is given by   	
\begin{align}
\Phi_d(2,0,r,s,\theta)&= \frac{1}{d^2}\Phi_d(1,0,r+1,s,\theta)\; \nonumber \\
&-\frac{2}{d}\Phi_d(1,1,r,s,\theta)+\frac{1}{d^2}(T_1(r,s,\theta)+T_2(r,s,\theta)),
\end{align}
where $B_{d_1,d_2}(r_1,s_1,r_2,s_2,\theta)$ is as defined in equation \eqref{bd1d2_simple}, 
\[T_1(r,s,\theta):=\sum\limits_{\substack{r_1+r_2=r \\ s_1+s_2=s \\ d_1+d_2=d \\ d_1,d_2>0}}\binom{r}{r_1}\binom{s}{s_1}d_1d_2^3 B_{d_1,d_2}(r_1,s_1,r_2,s_2,\theta),\]
\[T_2(r,s,\theta):=\sum\limits_{\substack{r_1+r_2=r \\ s_1+s_2=s \\ d_1+d_2=d \\ d_1,d_2>0}}\binom{r}{r_1}\binom{s}{s_1}d_1d_2^3 \widetilde{B}_{d_1,d_2}(r_1,s_1,r_2,s_2,\theta),\]
and
\[\widetilde{B}_{d_1,d_2}(r_1,s_1,r_2,s_2,\theta):=\sum_{i=0}^3 \Phi_{d_1}(1,0,r_1,s_1,i)  \times N_{d_2}^{\mathbb{P}^3, \mathrm{Planar}}(r_2,s_2,\theta+3-i).\]
\end{lemma}

\begin{rem}
The number $\Phi_d(1,j,r,s,\theta)$ for $j >1$ and $\Phi_d(2,j,r,s,\theta)$ for $j>0$ can be computed without any further effort; we have 
not presented the formulas since they are not needed for the Euler class computation. 
\end{rem}

Before we start proving these Lemmas, let us first recall an important result about $c_1(\mathbb{L^*})$. 

\begin{lemma}
\label{c1_divisor_ionel}
On $\overline{\mathcal{M}}_{0,1}^{\textnormal{Planar}}(\mathbb{P}^3, d)$, the following equality of divisors holds:
\begin{align*}
c_1(\mathbb{L}^*) & = \frac{\mathcal{H}_L}{d^2} -\frac{2}{d} \textnormal{ev}^*(H) +  \frac{1}{d^2}\sum_{\substack{d_1 + d_2 = d, \\ d_1, d_2 \neq 0}} d_2^2\mathcal{B}_{d_1, d_2},
\end{align*}
where $\mathcal{B}_{d_1, d_2}$ denotes the
boundary stratum corresponding to a bubble map of 
degree $d_1$ curve and degree $d_2$ curve with the marked point
lying on the degree $d_1$ component. 
\end{lemma}
\noindent \textbf{Proof:} This lemma is proved in (\cite{ionel}, Lemma 2.3) for  $\overline{\mathcal{M}}_{0,1}(\mathbb{P}^3, d)$. 
The corresponding statement for $\overline{\mathcal{M}}_{0,1}^{\textnormal{Planar}}(\mathbb{P}^3, d)$ 
follows immediately by pulling pack the relationship via the natural map $\pi_F$ from 
$\overline{\mathcal{M}}_{0,1}^{\textnormal{Planar}}(\mathbb{P}^3, d)$ to 
$\overline{\mathcal{M}}_{0,1}(\mathbb{P}^3, d)$ (the map that forgets the plane). \qed \\
\hf \hf We are now ready to prove the Lemmas that involve the computation of $\Phi_d(i,j,r,s,\theta)$. \\

\noindent \textbf{Proofs of Lemmas \ref{l0}, \ref{l1} and \ref{l2}:} Let us start by proving Lemma \ref{l0}. 
We recall that 
\begin{align*}
\Phi_d(0,j,r,s, \theta) &= 
\langle (\textnormal{ev}^*H)^j,\; \overline{\mathcal{M}}_{0,1}^{\mathrm{Planar}}(\mathbb{P}^3,d)\cap 
\mathcal{H}_L^r \cap \mathcal{H}_p^s\cap a^{\theta}\rangle.
\end{align*}
This number is zero unless $r+2s+\theta + j = 3d+3$. Let us start by considering the case when $j=0$. 
Let us assume $r+2s+\theta = 3d+3$ (otherwise the number is zero). 
In that case, we note that 
\begin{align*}
\Phi_d(0,0,r,s, \theta) &= 
\langle (\textnormal{ev}^*H)^0,\; \overline{\mathcal{M}}_{0,1}^{\mathrm{Planar}}(\mathbb{P}^3,d)\cap 
\mathcal{H}_L^r \cap \mathcal{H}_p^s\cap a^{\theta}\rangle \\ 
& =\langle \mathbf{1},\; \overline{\mathcal{M}}_{0,1}^{\mathrm{Planar}}(\mathbb{P}^3,d)\cap 
\mathcal{H}_L^r \cap \mathcal{H}_p^s\cap a^{\theta} \rangle \\ 
& = 0. 
\end{align*}
The last equality holds because 
the intersection is occurring inside $\overline{\mathcal{M}}_{0,0}^{\mathrm{Planar}}(\mathbb{P}^3,d)$ and $r+2s+\theta = 3d+3$.\\ 
\hf \hf Next, let us consider the case when $j =1$. 
Let us assume $r+2s+\theta+1 = 3d+3$ (otherwise the number is zero).
Let us consider the forgetful map  
\begin{align*}
\delta&: \overline{\mathcal{M}}_{0,1}^{\mathrm{Planar}}(\mathbb{P}^3,d) \cap \textnormal{ev}^*(H) \longrightarrow \overline{\mathcal{M}}_{0,0}^{\mathrm{Planar}}(\mathbb{P}^3,d).  
\end{align*}
We note that the degree of this map is $d$ (since a degree $d$ curve intersects a plane at $d$ points). Hence 
\begin{align*}
\Phi_d(0,1,r,s, \theta) &= 
\langle (\textnormal{ev}^*H),\; \overline{\mathcal{M}}_{0,1}^{\mathrm{Planar}}(\mathbb{P}^3,d)\cap 
\mathcal{H}_L^r \cap \mathcal{H}_p^s\cap a^{\theta}\rangle \\ 
& = \langle \mathbf{1},\; \Big(\overline{\mathcal{M}}_{0,1}^{\mathrm{Planar}}(\mathbb{P}^3,d) \cap \textnormal{ev}^*(H) \Big) \cap 
\mathcal{H}_L^r \cap \mathcal{H}_p^s\cap a^{\theta}\rangle \\
& = \textnormal{deg}(\delta)  |\overline{\mathcal{M}}_{0,0}^{\mathrm{Planar}}(\mathbb{P}^3,d)\cap 
\mathcal{H}_L^r \cap \mathcal{H}_p^s\cap a^{\theta}| \\ 
& = dN_d^{\mathbb{P}^3, \mathrm{Planar}}(r,s,\theta).
\end{align*}
\hf \hf Next, let us consider the case when $j =2$. 
Let us assume $r+2s+\theta+2 = 3d+3$ (otherwise the number is zero).
Let us consider the forgetful map  
\begin{align*}
\delta&: \overline{\mathcal{M}}_{0,1}^{\mathrm{Planar}}(\mathbb{P}^3,d) \cap \textnormal{ev}^*(H^2) \longrightarrow 
\overline{\mathcal{M}}_{0,0}^{\mathrm{Planar}}(\mathbb{P}^3,d) \cap \mathcal{H}_L.  
\end{align*}
We note that the degree of this map is one since in the first case we are considering a curve and a marked point that 
goes to a line ($H^2$), while in the second case we are considering a curve whose image intersects a line. Hence 
\begin{align*}
\Phi_d(0,2,r,s, \theta) &= 
\langle (\textnormal{ev}^*H^2),\; \overline{\mathcal{M}}_{0,1}^{\mathrm{Planar}}(\mathbb{P}^3,d)\cap 
\mathcal{H}_L^r \cap \mathcal{H}_p^s\cap a^{\theta}\rangle \\ 
& = \langle \mathbf{1},\; \Big(\overline{\mathcal{M}}_{0,1}^{\mathrm{Planar}}(\mathbb{P}^3,d) \cap \textnormal{ev}^*(H^2) \Big) \cap 
\mathcal{H}_L^r \cap \mathcal{H}_p^s\cap a^{\theta}\rangle \\
& = \textnormal{deg}(\delta)  |\overline{\mathcal{M}}_{0,0}^{\mathrm{Planar}}(\mathbb{P}^3,d)\cap 
\mathcal{H}_L^{r+1} \cap \mathcal{H}_p^s\cap a^{\theta}| \\ 
& = N_d^{\mathbb{P}^3, \mathrm{Planar}}(r+1,s,\theta).
\end{align*}
\hf \hf Finally, let us consider the case when $j =3$. 
Let us assume $r+2s+\theta+3 = 3d+3$ (otherwise the number is zero).
Let us consider the forgetful map  
\begin{align*}
\delta&: \overline{\mathcal{M}}_{0,1}^{\mathrm{Planar}}(\mathbb{P}^3,d) \cap \textnormal{ev}^*(H^3) \longrightarrow 
\overline{\mathcal{M}}_{0,0}^{\mathrm{Planar}}(\mathbb{P}^3,d) \cap \mathcal{H}_p.  
\end{align*}
We note that the degree of this map is one since in the first case we are considering a curve and a marked point that 
goes to a point ($H^3$), while in the second case we are considering a curve whose image passes through a point. Hence 
\begin{align*}
\Phi_d(0,3,r,s, \theta) &= 
\langle (\textnormal{ev}^*H^3),\; \overline{\mathcal{M}}_{0,1}^{\mathrm{Planar}}(\mathbb{P}^3,d)\cap 
\mathcal{H}_L^r \cap \mathcal{H}_p^s\cap a^{\theta}\rangle \\ 
& = \langle \mathbf{1},\; \Big(\overline{\mathcal{M}}_{0,1}^{\mathrm{Planar}}(\mathbb{P}^3,d) \cap \textnormal{ev}^*(H^3) \Big) \cap 
\mathcal{H}_L^r \cap \mathcal{H}_p^s\cap a^{\theta}\rangle \\
& = \textnormal{deg}(\delta)  |\overline{\mathcal{M}}_{0,0}^{\mathrm{Planar}}(\mathbb{P}^3,d)\cap 
\mathcal{H}_L^{r} \cap \mathcal{H}_p^{s+1}\cap a^{\theta}| \\ 
& = N_d^{\mathbb{P}^3, \mathrm{Planar}}(r,s+1,\theta).
\end{align*}
For $j>3$, let us assume $r+2s+\theta+j=3d+3$, since $ev^* H^j=0$ for all $j>3$. Therefore, we have $\Phi_d(0,j,r,s, \theta)=0$.

\vspace{.2cm}
\hf \hf Let us now prove Lemma \ref{l1}. 
We recall that 
\begin{align*}
\Phi_d(1,j,r,s, \theta) &= 
\langle c_1(\mathbb{L}^*) (\textnormal{ev}^*H)^j,\; \overline{\mathcal{M}}_{0,1}^{\mathrm{Planar}}(\mathbb{P}^3,d)\cap 
\mathcal{H}_L^r \cap \mathcal{H}_p^s\cap a^{\theta}\rangle.
\end{align*}
This number is zero unless $r+2s+\theta + 1+j = 3d+3$. Let us start by considering the case when $j=0$. 
Let us assume $r+2s+\theta +1 = 3d+3$ (otherwise the number is zero). We note that 
\begin{align}
\langle \mathcal{H}_L,\; \overline{\mathcal{M}}_{0,1}^{\mathrm{Planar}}(\mathbb{P}^3,d)\cap 
\mathcal{H}_L^r \cap \mathcal{H}_p^s\cap a^{\theta}\rangle & = 0, \label{itc1} \\ 
\langle \textnormal{ev}^*(H),\; \overline{\mathcal{M}}_{0,1}^{\mathrm{Planar}}(\mathbb{P}^3,d)\cap 
\mathcal{H}_L^r \cap \mathcal{H}_p^s\cap a^{\theta}\rangle & = \Phi_d(0, 1, r, s, \theta) \nonumber \\ 
 & = d N_d^{\mathbb{P}^3, \mathrm{Planar}}(r, s, \theta).\label{itc2} 
\end{align}
\hf \hf Next, let us consider the boundary divisor  
$\mathcal{B}_{d_1, d_2}$, the
boundary stratum of $\overline{\mathcal{M}}_{0,1}^{\mathrm{Planar}}(\mathbb{P}^3,d)$ corresponding to a bubble map of 
degree $d_1$ curve and degree $d_2$ curve, with the marked point
lying on the degree $d_1$ component. Let us now compute the degree of the divisor 
\begin{align}
B_{d_1, d_2}(r, s, \theta)&:= \textnormal{deg}\Big(\mathcal{B}_{d_1, d_2} \cap \textnormal{ev}^*(H^j) \cap \mathcal{H}_L^r \cap \mathcal{H}_p^s \cap a^{\theta}\Big). 
\label{bd1d2_deg}
\end{align}
As per the convention decided in equation \eqref{convention_deg}, we formally declare the degree to be zero unless 
$r+2s+\theta+1+j = 3d+3$. First, let us compute the number of (ordered) two component rational curves of type $(d_1, d_2)$ that 
lies inside $\overline{\mathcal{M}}_{0,0}^{\mathrm{Planar}}(\mathbb{P}^3,d) \cap a^{\theta}$. 
This is computed in (\cite{ritwik}), in the proof of Theorem 3.3 (and in equation (2.2)), given by 
\begin{align}
B_{d_1, d_2}(r_1, s_1, r_2, s_2, \theta) & = 
\sum_{i=0}^3 N_{d_1}^{\mathbb{P}^3, \mathrm{Planar}}(r_1,s_1,i)\times N_{d_2}^{\mathbb{P}^3, \mathrm{Planar}}(r_2,s_2,\theta+3-i). \label{bd1d2_simple}
\end{align}
We note there that each element of $B_{d_1, d_2}(r_1, s_1, r_2, s_2, \theta)$ corresponds to $d_1 d_2$ bubble  maps in 
$\overline{\mathcal{M}}_{0,0}^{\mathrm{Planar}}(\mathbb{P}^3,d)$, because there are $d_1 d_2$ choices for the nodal point of the domain. 
Furthermore, each such bubble map corresponds to $d_1$ bubble maps in $\overline{\mathcal{M}}_{0,1}^{\mathrm{Planar}}(\mathbb{P}^3,d) \cap \textnormal{ev}^*(H)$ 
since each degree $d_1$ curve intersects a hyperplane in $d_1$ points. 
Hence, 
\begin{align}
B_{d_1, d_2}(r, s, \theta)&= \sum_{\substack{r_1  + r_2 = r, \\ s_1 + s_2 = s}} \binom{r}{r_1}\binom{s}{s_1} d_1^2 d_2B_{d_1, d_2}(r_1, s_1, r_2, s_2, \theta).
\label{bd1d2_full}
\end{align}
By equation \eqref{bd1d2_simple} and \eqref{bd1d2_full}, we conclude that $B_{d_1, d_2}(r, s, \theta)$ is zero if $1+r+2s+\theta = 3d+3$. 
Hence, equations \eqref{itc1}, \eqref{itc2}, \eqref{bd1d2_full} and Lemma \ref{c1_divisor_ionel} imply the result of Lemma \ref{l1} for $j=0$. \\ 
\hf \hf Next, let us prove Lemma \ref{l1} for the case when $j=1$. 
Let us assume $r+2s+\theta +1+1 = 3d+3$ (otherwise the number is zero). We note that 
\begin{align}
\langle \mathcal{H}_L \cdot \textnormal{ev}^*(H),\; \overline{\mathcal{M}}_{0,1}^{\mathrm{Planar}}(\mathbb{P}^3,d)\cap 
\mathcal{H}_L^r \cap \mathcal{H}_p^s\cap a^{\theta}\rangle & = \Phi_d(0,1,r+1, s, \theta), \label{itc12} \\ 
\langle \textnormal{ev}^*(H^2),\; \overline{\mathcal{M}}_{0,1}^{\mathrm{Planar}}(\mathbb{P}^3,d)\cap 
\mathcal{H}_L^r \cap \mathcal{H}_p^s\cap a^{\theta}\rangle & = \Phi_d(0, 2, r, s, \theta). \label{itc22} 
\end{align}
 We note that equation \eqref{bd1d2_full} is true irrespective of the value of $r$ and $s$. Hence, using 
equations \eqref{itc12}, \eqref{itc22}, \eqref{bd1d2_full} and Lemma \ref{c1_divisor_ionel}, we obtain the result of 
Lemma \ref{l1} for the case when $j=1$. 

\vspace{.2cm}
\hf \hf Finally, let us prove Lemma \ref{l2}. 
We recall that 
\begin{align*}
\Phi_d(2,j,r,s, \theta) &= 
\langle c_1(\mathbb{L}^*)^2 (\textnormal{ev}^*H)^j,\; \overline{\mathcal{M}}_{0,1}^{\mathrm{Planar}}(\mathbb{P}^3,d)\cap 
\mathcal{H}_L^r \cap \mathcal{H}_p^s\cap a^{\theta}\rangle.
\end{align*}
This number is zero unless $r+2s+\theta + 2+j = 3d+3$. We will only consider the case when $j=0$. 
Let us assume $r+2s+\theta +2 = 3d+3$ (otherwise the number is zero). We note that 
\begin{align}
\langle c_1(\mathbb{L}^*)\cdot \mathcal{H}_L,\; \overline{\mathcal{M}}_{0,1}^{\mathrm{Planar}}(\mathbb{P}^3,d)\cap 
\mathcal{H}_L^r \cap \mathcal{H}_p^s\cap a^{\theta}\rangle & = \Phi_d(1,0,r+1,s, \theta), \label{itc13} \\ 
\langle c_1(\mathbb{L}^*)\cdot \textnormal{ev}^*(H),\; \overline{\mathcal{M}}_{0,1}^{\mathrm{Planar}}(\mathbb{P}^3,d)\cap 
\mathcal{H}_L^r \cap \mathcal{H}_p^s\cap a^{\theta}\rangle & = \Phi_d(1, 1, r, s, \theta). \label{itc23} 
\end{align}
Next we will show that 
\begin{align}
\sum_{\substack{d_1+d_2 = d, \\ d_1, d_2 \neq 0}} d_2^2 \Big\langle c_1(\mathbb{L}^*)\cdot \mathcal{B}_{d_1, d_2},\; 
\overline{\mathcal{M}}_{0,1}^{\mathrm{Planar}}(\mathbb{P}^3,d)\cap 
\mathcal{H}_L^r \cap \mathcal{H}_p^s\cap a^{\theta}\Big\rangle &= T_1(r,s,\theta) + T_2(r, s, \theta). \label{itc_b}
\end{align}
We note that equations \eqref{itc13}, \eqref{itc23}, \eqref{itc_b} and Lemma \ref{c1_divisor_ionel} implies Lemma \ref{l2}. \\ 
\hf \hf Let us now prove equation \eqref{itc_b}. Let us consider the map  
\begin{align*}
\pi&:\Big(\mathcal{B}_{d_1, d_2} \cap \mathcal{H}_L^{r_1} \cap \mathcal{H}_p^{s_1}\Big) \longrightarrow  
\overline{\mathcal{M}}_{0,1}^{\mathrm{Planar}}(\mathbb{P}^3,d_1)\cap 
\mathcal{H}_L^{r_1} \cap \mathcal{H}_p^{s_1},
\end{align*}
that maps to the degree $d_1$ component. Let $\mathbb{L}_1$ denotes the pullback of the universal tangent bundle 
over $\overline{\mathcal{M}}_{0,1}^{\mathrm{Planar}}(\mathbb{P}^3,d_1)$ to $\mathcal{B}_{d_1, d_2}$ via the map $\pi$. 
By (\cite{ionel}) (equation (2.10), Page 29), we conclude that on 
$\Big(\mathcal{B}_{d_1, d_2} \cap \mathcal{H}_L^{r_1} \cap \mathcal{H}_p^{s_1}\Big)$, we have 
the equality of divisors \begin{align*}
c_1(\mathbb{L}^*)|_{\mathcal{B}_{d_1, d_2}}&= c_1(\mathbb{L}_1^*) + \mathcal{G}. 
\end{align*}
Now let us compute each of the terms. Let us now consider the space 
\begin{align*}
\overline{\mathcal{M}}_{0,0}^{\mathrm{Planar}}(\mathbb{P}^3,d_2) \cap \mathcal{H}_L^{r_2} \cap \mathcal{H}_p^{s_2} \cap a^{\theta}. 
\end{align*}
Hence, we conclude that 
\begin{align*}
& \Big\langle c_1(\mathbb{L}_1^*),\; 
\overline{\mathcal{M}}_{0,1}^{\mathrm{Planar}}(\mathbb{P}^3,d)\cap 
\mathcal{H}_L^r \cap \mathcal{H}_p^s \cap a^{\theta} \Big\rangle \\ 
& = \sum_{\substack{r_1+r_2 = r, \\s_1+s_2 = s }} \binom{r}{r_1} \binom{s}{s_1} d_1 d_2\Big \langle  c_1(\mathbb{L}_{1}^*) \cdot \Delta_{\widehat{\mathbb{P}}^3},\\ 
& \qquad 
\overline{\mathcal{M}}_{0,1}^{\mathrm{Planar}}(\mathbb{P}^3,d_1)\cap \mathcal{H}_L^{r_1} \cap \mathcal{H}_p^{s_1} \times 
\overline{\mathcal{M}}_{0,0}^{\mathrm{Planar}}(\mathbb{P}^3,d_2)\cap \mathcal{H}_L^{r_2} \cap \mathcal{H}_p^{s_2} \cap a^{\theta} \Big\rangle \\ 
& = \sum_{\substack{r_1+r_2 = r, \\s_1+s_2 = s }} \binom{r}{r_1} \binom{s}{s_1} d_1 d_2 \widetilde{B}_{d_1, d_2}(r_1, s_1, r_2, s_2, \theta).
\end{align*}
Here $\mathbb{L}_{1}^*$ denotes the universal cotangent bundle over the first marked point (over the moduli space of degree $d_1$ curves) 
and $\Delta_{\widehat{\mathbb{P}}^3}$denotes the diagonal of $\widehat{\mathbb{P}}^3 \times \widehat{\mathbb{P}}^3$. 
We note that the class of the diagonal is given by $\sum_{i=0}^3\pi_1^*a^i\pi_2^*a^{3-i}$, 
where $\pi_{i}:\widehat{\mathbb{P}}^2 \times \widehat{\mathbb{P}}^2 \longrightarrow \widehat{\mathbb{P}}^2$ is the 
projection map to the $i^{\textnormal{th}}$ factor. \\  
\hf \hf Similarly, we note that
\begin{align*}
& \Big\langle \mathcal{G},\; 
\overline{\mathcal{M}}_{0,1}^{\mathrm{Planar}}(\mathbb{P}^3,d)\cap 
\mathcal{H}_L^r \cap \mathcal{H}_p^s \cap a^{\theta} \Big\rangle \\ 
& = \sum_{\substack{r_1+r_2 = r, \\s_1+s_2 = s }} \binom{r}{r_1} \binom{s}{s_1} d_1 d_2\Big \langle  
\Delta_{\widehat{\mathbb{P}}^3},\\ 
& \qquad 
\overline{\mathcal{M}}_{0,0}^{\mathrm{Planar}}(\mathbb{P}^3,d_1)\cap \mathcal{H}_L^{r_1} \cap \mathcal{H}_p^{s_1} \times 
\overline{\mathcal{M}}_{0,0}^{\mathrm{Planar}}(\mathbb{P}^3,d_2)\cap \mathcal{H}_L^{r_2} \cap \mathcal{H}_p^{s_2} \cap a^{\theta}\Big\rangle \\ 
& = \sum_{\substack{r_1+r_2 = r, \\s_1+s_2 = s }} d_1 d_2 \binom{r}{r_1} \binom{s}{s_1} B_{d_1, d_2}(r_1, s_1, r_2, s_2, \theta).
\end{align*}
This proves the claim. \qed


\section{Transversality and degenerate contribution to the Euler class} 
\label{Degenerate_Euler}
Let us start by showing that the section of the bundle, considered in equation 
\eqref{du_transverse} is transverse to zero (restricted to $\mathcal{M}$). 
This follows from the fact that restricted to each fibre, the section is transverse to zero. 
Fibre wise transversality is proved in 
\cite[Lemma 8.2,  Page 105]{BSRV}
and it is also used in \cite{g2p2and3}. We  will now justify the multiplicity.\\ 
\hf \hf Let us recapitulate the notation of section \ref{notation}.  
We have defined 
\begin{align*}
\pi: \overline{\mathcal{M}}_{0,1}^{\mathrm{Planar}}(\mathbb{P}^3,d) & \longrightarrow \widehat{\mathbb{P}}^3 
\end{align*}
to be the fiber bundle, such that the fiber over each point is given by 
\begin{align*}
\pi^{-1}([\eta]) & := \overline{\mathcal{M}}_{0,1}(\mathbb{P}^2_{\eta}, d). 
\end{align*}
Let us abbreviate the Euler class of the rank two bundle considered in section \ref{Euler_class_computation} 
as $E$, i.e. 
\begin{align*}
E&:= e(\mathbb{L}^* \otimes \textnormal{ev}^*W). 
\end{align*}
Note that this is a (complex) degree two cohomology class (codimension two cycle) in $\overline{\mathcal{M}}_{0,1}^{\mathrm{Planar}}(\mathbb{P}^3,d)$; 
it is not a number. \\ 
\hf \hf Next, let us denote $B \subset \overline{\mathcal{M}}_{0,1}^{\mathrm{Planar}}(\mathbb{P}^3,d)$ to be the boundary class; i.e. 
it is the closure of the space of bubble maps of type $(d_1, d_2)$ such that the marked point is the nodal point 
(or more precisely, it is a map from a wedge of three spheres, where the middle component is constant and the marked point lies on it 
while the other two components are of degree $d_1$ and $d_2$). \\ 
\hf \hf Finally, let $C \subset \overline{\mathcal{M}}_{0,1}^{\mathrm{Planar}}(\mathbb{P}^3,d)$ denote the class determined by the 
closure of the space of cuspidal curves (the cusp being  on the marked point). Let 
\[ \psi:\overline{\mathcal{M}}_{0,1}^{\mathrm{Planar}}(\mathbb{P}^3,d) \longrightarrow \mathbb{L}^* \otimes \textnormal{ev}^*W \] 
denote the section that corresponds to taking the derivative at the marked point 
(i.e. the section defined in equation \eqref{du_transverse}). We note that the section vanishes on $C$ and $B$. 
Hence, set theoretically, $\psi^{-1}(0)$ is the union of $C$ and $B$. Hence, as cycles, 
we conclude that 
\begin{align}
E & = m_1 [C] + m_2 [B], \label{E_cycle}  
\end{align}
where $m_1$ and $m_2$ are integers. Since the section $\psi$ vanishes transversally on the cuspidal curves, 
we conclude that $m_1 =1$. 
Let us now show that $m_2=1$. Consider the cycle  
\[ [\mathcal{Z}]:= \mathcal{H}_l^{3d-2} a^3. \] 
Now intersect the left hand side and right hand side of equation \eqref{E_cycle} with $\mathcal{Z}$ to conclude that 
\begin{align*}
E\cdot [\mathcal{Z}] &= [C]\cdot [\mathcal{Z}] + m_2 [B]\cdot [\mathcal{Z}]. 
\end{align*}
We now note that intersecting with $a^3$ 
corresponds to fixing a $\mathbb{P}^2$. Hence we are in the situation considered in \cite{BSRV}, where 
we show that  
\begin{align*}
E\cdot [\mathcal{Z}] &= [C]\cdot [\mathcal{Z}] + 1 \times [B]\cdot [\mathcal{Z}]. 
\end{align*}
Hence, $m_2 =1$. 
Let us now consider the following cycle in $\overline{\mathcal{M}}_{0,1}^{\mathrm{Planar}}(\mathbb{P}^3,d)$
\begin{align*}
[\mathcal{Z}] &:= \mathcal{H}_L^r\cdot \mathcal{H}_p^s \cdot a^{\theta}. 
\end{align*}
Choose $r,s$ and $\theta$ such that dimension of $\mathcal{Z}$ is two (i.e. $r+2s+\theta = 3d+1$). 
Intersecting the left hand side and right hand side of equation \eqref{E_cycle} with $\mathcal{Z}$,  
we get the following equality of numbers 
\begin{align*}
E\cdot [\mathcal{Z}] & = [C]\cdot [\mathcal{Z}] + [B] \cdot [\mathcal{Z}]. 
\end{align*}
This is precisely equation \eqref{Euler_Cd}. \\

\section{Low degree checks}
\label{low_deg_chck}
In this section we subject our formula to certain low degree checks. All these numbers have been computed using 
our mathematica program. We will abbreviate $C_d^{\mathbb{P}^3, \mathrm{Planar}}(r,s)$ by $C_d(r,s)$.\\
\hf \hf First of all our formula gives us the value of zero for $C_d(r,s)$ when $d =2$. 
This is as geometrically as expected since there are no conics with a cusp. \\ 
\hf \hf Next, in (\cite{ND}), N. Das and the first author compute the following numbers: 
what is $N_d(A_1^{\delta} A_2, r, s)$, the number 
of planar degree $d$ curves in $\mathbb{P}^3$, passing through $r$ lines and 
$s$ points,
that have $\delta$ (ordered) nodes and one cusp, for all $\delta \leq 2$. 
Note that here $r +  2s = \delta + 2$.
For 
$d =3$, and $\delta =0$, this number should be the same as the characteristic number of genus zero planar cubics in $\mathbb{P}^3$ with a cusp, 
i.e. $C_d(r, s)$.  We have verified that is indeed the case. We tabulate the numbers for the readers convenience: 
\begin{align*}
C_3(10, 0)& = 17760, \qquad C_3(8, 1)  = 2064,  \\ 
 C_3(6, 2)& = 240 \qquad \textnormal{and} \qquad    
C_3(4, 3) = 24.
\end{align*}
These numbers are the same as $N_d(A_1^{\delta} A_2, r, s)$ for $d=3$ and $\delta =0$. \\ 
\hf \hf Next, we note that when $d=4$ and $\delta =2$, the number $\frac{1}{\delta !}N_d(A_1^{\delta} A_2, r, s)$ 
is same as the characteristic number of genus zero planar quartics in $\mathbb{P}^3$ with a cusp, i.e.   
$C_d(r, s)$.  We have verified that fact. The numbers are  
\begin{align*}
C_4(13, 0)& = 10613184, \qquad C_4(11, 1) = 760368, \\ 
C_4(9, 2) & = 49152 \qquad \textnormal{and} \qquad  
C_4(7, 3) = 2304.
\end{align*}
These numbers are the same as $\frac{1}{2!}N_d(A_1^{\delta} A_2, r, s)$ for $d=4$ and $\delta =2$. 
We have to divide out by a factor of $\delta!$ because in the definition of 
$N_d(A_1^{\delta} A_2, r, s)$, the nodes are ordered.

\section{Acknowledgement}
The first author would like to acknowledge the External Grant 
he has obtained, namely 
MATRICS (File number: MTR/2017/000439) that has been sanctioned by the Science and Research Board (SERB).

\end{document}